\documentclass[12pt]{amsart}
\newtheorem{theorem}{Theorem}

\theoremstyle{definition}

\theoremstyle{remark}

\newcommand{\NN}{ {\mathbb N} }

\newcommand{\CC}{{\mathbb C}}

\newcommand{\ff}{\varphi}

\newcommand{\cB}{\mathcal{B}}
\newcommand{\cA}{\mathcal{A}}

\newcommand{\id}{\text{id}}

\newcommand{\cC}{\mathcal{C}}

\newcommand{\tovert}{\to}

\newcommand{\xhk}{x^{(k)}}

\newcommand{\Shk}{S^{(k)}}
\newcommand{\hi}{^{[i]}}

\renewcommand{\Im}{\text{{\rm Im\,}}}

\begin{document}

\title[Berry Esseen for multivariate free CLT]
{On the rate of convergence and Berry-Esseen type theorems
for a multivariate free central limit theorem}

\author[R. Speicher]{{Roland Speicher}
$^{(\dagger)}$}
\thanks{$^\dagger\,$Research supported by Discovery and LSI grants from NSERC (Canada) and by
a Killam Fellowship from the Canada Council for the Arts}
\address{Queen's University, Department of Mathematics and Statistics,
Jeffery Hall, Kingston, ON, K7L 3N6, Canada} \email{speicher@mast.queensu.ca}

\begin{abstract}
We address the question of a Berry Esseen type theorem for the speed of convergence in a
multivariate free central limit theorem. For this, we estimate the difference between the
operator-valued Cauchy transforms of the normalized partial sums in an operator-valued
free central limit theorem and the Cauchy transform of the limiting operator-valued
semicircular element.
\end{abstract}

\thanks{This project was initiated by disussions with Friedrich G\"oetze during
my visit at the University of Bielefeld in November 2006. I thank the Department of
Mathematics and in particular the SFB 701 for its generous hospitality and Friedrich
G\"otze for the invitation and many interesting discussions.\\
I also thank Uffe Haagerup for pointing out how ideas from \cite{HST} can be used to
improve the results from an earlier version of this paper.} \maketitle

\section{Introduction}
The free central limit theorem (due to Voiculescu \cite{VoiZGW} in the one-dimensional
case, and to Speicher \cite{SpeZGW} in the multivariate case) is one of the basic results
in free probability theory. Investigations on the speed of convergence to the limiting
semicircular distribution, however, were taken up only recently. In the classical
context, the analogous question is answered by the famous Berry-Esseen theorem, which
states, in its simplest version, the following: If $X_i$ are i.i.d. random variables,
with mean zero and variance 1, then the distance between $S_n:=(X_1+\cdots+X_n)/\sqrt n$
and a normal variable $\gamma$ of mean zero and variance 1 can be estimated in terms of
the Kolmogorov distance $\Delta$ by
$$\Delta(S_n,\gamma)\leq C\frac 1{\sqrt n} \rho,$$
where $C$ is a constant and $\rho$ is the absolute third moment of the variables $x_i$.

The question for a free analogue of the Berry-Esseen estimate in the case of one random
variable was answered by Chistyakov and G\"otze \cite{CG}: If $x_i$ are free identically
distributed random variables with mean zero and variance 1, then the distance between
$S_n:=(X_1+\cdots+X_n)/\sqrt n$ and a semicircular variable $s$ of mean zero and variance
1 can, under the assumption of finite fourth moment, be estimated as
$$\Delta(S_n,s)\leq c \frac {\vert m_3\vert +\sqrt {m_4}}{\sqrt n},$$
where $c>0$ is an absolute constant, and $m_3$ and $m_4$ are the third and fourth moment,
respectively, of the $x_i$. (Independently, the same kind of question was considered,
under the more restrictive assumption of compact support for the $x_i$, by Kargin
\cite{K}.)

In this paper we want to address the multivariate version of a free Berry-Esseen theorem.
In contrast to the classical situation, the multivariate situation is of a quite
different nature than the one-dimensional case, because we have to deal with
non-commuting operators and all the analytical tools, which are available in the
one-dimensional case, break down. However, we are able to deal with this situation by
invoking recent ideas of Haagerup and Thorbjornsen \cite{HT,HST}, in particular, their
linearization trick which allows to reduce the multivariate (scalar-valued) to an
analogous one-dimensional operator-valued problem. Estimates for the operator-valued
Cauchy transform of this operator-valued operator are quite similar to estimates in the
scalar-valued case. Actually, on the level of deriving equations for these Cauchy
transforms we can follow ideas which are used for dealing with speed of convergence
questions for random matrices; here we are inspired in particular by the work of G\"otze
and Tikhomirov \cite{GT}, but see also \cite{B1,B2}. Our main theorem on the speed of
convergence in an operator-valued free central limit theorem is the following.

\begin{theorem}\label{thm:operator-valued}
Let $1\in\cB\subset\cA$, $E:\cA\to\cB$ be an operator-valued probability space. Consider
selfadjoint $X_1,X_2,\dots \in\cA$ which are free with respect to $E$ and have identical
$\cB$-valued distribution. Assume that the first moments vanish,
$$E[X_i]=0$$
and let
$$\eta:\cB\to\cB,\qquad \eta(b)=E[X_ibX_i]$$
be their covariance. Denote
$$\alpha_2:=\sup_{b\in\cB\atop \Vert b\Vert =1} \Vert E[X_ibX_i]\Vert=\Vert \eta\Vert$$
and
$$\alpha_4:=\sup_{b\in\cB\atop \Vert b\Vert=1}\Vert E[X_ibX_iX_ib^*X_i]\Vert.$$
Consider now the normalized sums
$$S_n:=\frac{X_1+\cdots+X_n}{\sqrt n}$$
and their $\cB$-valued Cauchy transforms
$$G_n(b):=E[\frac 1{b-S_n}] \qquad (b\in\cB_+)$$
on the ``upper half plane'' $\cB_+$ in $\cB$,
$$\cB_+:=\{b\in \cB\mid \Im b\geq 0 \text{ and $\Im b$ invertible}\}.$$
By $G$ we denote the operator-valued Cauchy transform of a $\cB$-valued semicircular
element with covariance $\eta$.

Then we have for all $b\in\cB_+$ and all $n\in\NN$ that
\begin{equation}
\Vert G_n(b)-G(b)\Vert\leq 4c_n(b)\left(\Vert b\Vert +\alpha_2 \cdot\Vert \frac 1{\Im
b}\Vert\right)\cdot \Vert \frac 1{\Im b}\Vert^2,
\end{equation}
where
$$c_n(b):=\frac 1{\sqrt n} \bigl\Vert\frac 1{\Im b}\bigr\Vert^3\sqrt{\alpha_2}
\cdot(2\alpha_2+\sqrt{\alpha_4+2\alpha_2^2})+\frac 1n \bigl\Vert\frac 1{\Im
b}\bigr\Vert^4\alpha_2^2.
$$
\end{theorem}

In the one-dimensional scalar case one can derive from such estimates corresponding
estimates for the Kolmogorov distance between the distribution of $S_n$ and the limiting
semicircle $s$. This relies on the fact that the Kolmogorov metric measures how close the
distribution functions of two measures are, and the Stieltjes inversion formula allows to
relate the distribution function with Cauchy transforms. (In the proof of the classical
Berry-Esseen theorem one follows a similar route, using Fourier transforms instead of
Cauchy transforms.) For the multivariate case, say of $d$ variables, where we would like
to say something about the speed of convergence of the $d$-tuple of partial sums
$(S_n^{(1)},\dots,S_n^{(d)})$ to the limiting semicircular family $(s_1,\dots,s_d)$,
there is no nice replacement for the distribution function, and we also do not know of a
canonical metric on joint distributions of several non-commuting variables which relates
directly with the above estimates for operator-valued Cauchy transforms.

However, there is a kind of replacement for this; namely, following again \cite{HST},
estimates for Cauchy transforms of linear combinations with operator-valued coefficients
of the variables $(S_n^{(1)},\dots,S_n^{(d)})$ should imply corresponding estimates for
any non-commutative scalar polynomial in those variables and from those one should be
able to estimate, for any selfadjoint non-commutative polynomial $p$, the Levy distance
between $p(S_n^{(1)},\dots,S_n^{(d)})$ and $p(s_1,\dots,s_d)$. However, one has to deal
with the following problem in such an approach: as is shown in \cite{HST} one can get the
Cauchy transform of a polynomial $p(s_1,\dots,s_d)$ as a corner of an operator-valued
Cauchy transform of a linear combination $P$, with matrix-valued coefficients, of
$s_1,\dots,s_d$; but, even if $p$ is a selfadjoint polynomial, the corresponding
matrix-valued operator $P$ is not selfadjoint, and thus our operator-valued estimates,
which were only shown for selfadjoint $X$, cannot be used directly for $P$; one would
have to reprove most of our statements also for $P$. It is conceivable that this can be
done in a similar manner as in \cite{HST}; as this approach is getting quite technical,
we will pursue the details in a forthcoming investigation.

Note that for proving such a kind of Berry-Esseen theorem for polynomials
$p(s_1,\dots,s_d)$ one also has to face another kind of question: estimates for the
difference of Cauchy transforms translate directly only in estimates for the Levy
distance between the corresponding measures; in order to get also estimates for the more
intuitive Kolmogorov distance one needs to know that the distribution of
$p(s_1,\dots,s_d)$ has a continuous density, in particular, has no atoms. We conjecture
that this is true for all non-commutative selfadjoint polynomials $p$ in a semicircular
family, but this seems to be a non-trivial problem. Note that the question of absence of
atoms can be seen as an analogue of the Zero-Divisor Theorem for the free group. We hope
to address this question in some future work.

The paper is organized as follows. In the next section we will first relate a
multivariate free central limit theorem with a one-dimensional operator-valued free
central limit theorem. The proof of Theorem \ref{thm:operator-valued} will be given in
Section \ref{sect:operator}.

\section{Multivariate free central limit theorem}
\subsection{Setting}\label{subsection:setting}
Let $\bigl(x^{(k)}_1\bigr)_{k=1}^d,\bigl(\xhk_2\bigr)_{k=1}^d,\dots$ be free and
identically distributed sets of $k$ selfadjoint random variables in some non-commutative
probability space $(\cC,\ff)$, such that the first moments vanish and the second moments
are given by a covariance matrix $\Sigma=(\sigma_{kl})_{k,l=1}^d$. We put
$$\Shk_n=\frac{\xhk_1+\dots+\xhk_n}{\sqrt n}.$$
We know \cite{SpeZGW} that $(S_n^{(1)},\dots,S_n^{(d)})$ converges in distribution for
$n\to\infty$ to a semicircular family $(s_1,\dots,s_d)$ of covariance $\Sigma$. We want
to analyze the rate of this convergence. We would like to get an estimate which involves
only small moments of the given variables. As we will see, the second and fourth moments
of our variables will show up in the estimates and we will use the upper bound
$$\beta_2:=\max_{k,l}\vert\sigma_{k,l}\vert=\max_{k,l}\ff(x_i^{(k)}x_i^{(l)})$$
for the second and the upper bound
$$\beta_4:=\max_{r,p,k,l}\vert\ff(x_i^{(r)}x_i^{(p)}x_i^{(k)}x_i^{(l)})\vert$$
for the fourth moments.

\subsection{Transition to operator-valued frame}\label{subsection:transition}
We will analyze the rate of convergence of the multivariate problem,
$$(S_n^{(1)},\dots,S_n^{(d)})\tovert (s_1,\dots,s_n)$$
by replacing this by an
one-dimensional operator-valued problem. The underlying idea for that is the
linearization trick \cite{HT,HST} that one can understand the joint distribution of
several scalar random variables by understanding the distribution of each operator-valued
linear combination of those random variables.

Let $\cB=M_N(\CC)$ and put $\cA:=M_N(\CC)\otimes \cC=M_N(\cC)$. Then $\cB\cong\cB\otimes
1\subset\cA$ is an operator-valued probability space with respect to the conditional
expectation
$$
E=\id\otimes \ff:\cB\otimes\cC\to\cB,\qquad b\otimes c\mapsto \ff(c)b.
$$

For some fixed $b_1,\dots,b_k\in M_N(\CC)$ we put
$$X_i:=\sum_{k=1}^d b_k\otimes x_i^{(k)}$$
and
$$S_n:=\sum_{k=1}^d b_k\otimes \Shk_n$$
Note that $X_1,X_2,\cdots$ are free with respect to $E$ and that we have
$$S_n=\frac{X_1+\cdots+X_n}{\sqrt n}.$$
The limit of $S_n$ is
$$s:=\sum_{k=1}^d b_k\otimes s_k,$$
which is an $\cB=M_N(\CC)$-valued semicircular element with covariance mapping
$\eta:\cB\to\cB$ given by
\begin{align*}
\eta(b)=E[s b\otimes 1 s]&=\sum_{k,l=1}^d E[b_k\otimes s_k\cdot  b\otimes 1\cdot  b_l
\otimes s_l]\\&= \sum_{k,l=1}^d b_k b b_l \ff(s_ks_l)= \sum_{k,l=1}^d b_k b b_l
\sigma_{kl}.
\end{align*}

We want to determine the rate of convergence for $S_n$ to $s$. We will do this in the
next section in the context of a general operator-valued free central limit theorem.

\section{Rate of convergence for operator-valued free central limit
theorem}\label{sect:operator}

\subsection{Setting} Let $1\in\cB\subset\cA$, $E:\cA\to\cB$ be an operator-valued probability
space. This means that $\cA$ is a von Neumann algebra, $\cB$ is a sub von Neumann
algebra, which contains the identity of $\cA$, and $E$ is a conditional expectation from
$\cA$ onto $\cB$, i.e., a linear map which satisfies the property
$$E[b_1 a b_2]=b_1E[a]b_2$$
for all $a\in\cA$ and $b_1,b_2\in\cB$.

Consider selfadjoint $X_1,X_2,\dots \in\cA$ which are free with respect to $E$ and have
identical $\cB$-valued distribution. Assume that the first moments vanish,
$$E[X_i]=0$$
and let
$$\eta:\cB\to\cB,\qquad \eta(b)=E[X_ibX_i]$$
be their covariance. We will need
$$\alpha_2:=\sup_{b\in\cB\atop \Vert b\Vert =1} \Vert E[X_ibX_i]\Vert=\Vert \eta\Vert$$
and
$$\alpha_4:=\sup_{b\in\cB\atop \Vert b\Vert=1}\Vert E[X_ibX_iX_ib^*X_i]\Vert.$$
Consider now the normalized sums
$$S_n:=\frac{X_1+\cdots+X_n}{\sqrt n}.$$
We know that $S_n$ converges in distribution to an operator-valued semicircular element $s$
with covariance $\eta$, see \cite{Spe}

We want to estimate the rate of this convergence. Let us denote by $\cB_+$ the ``upper
half plane'' in $\cB$, i.e.,
$$\cB_+:=\{b\in \cB\mid \Im b\geq 0 \text{ and $\Im b$ invertible}\}.$$
We consider, for $b\in\cB_+$, the resolvents
$$R_n(b):=\frac 1{b-S_n},\qquad R(b):=\frac 1{b-s}$$
and the Cauchy transforms
$$G_n(b):=E[R_n(b)],\qquad G(b):=E[R(b)].$$
$G_n$ and $G$ are analytic functions in $\cB_+$.

\subsection{The main estimates}
We will show that $G_n(b)$ converges to $G(b)$, where we have good control over the
difference in terms of $n$ and $b$. The idea for showing this is the same as in
\cite{HT}. First we show that $G_n$ satisfies an approximate version of an equation
satisfied by $G$ and then we show that this actually implies that $G_n$ and $G$ must be
close to each other.

Let us start with deriving the equations for $G$ and $G_n$.

Since $s$ is an operator-valued semicircular element with covariance $\eta$ we know
\cite{Voi,Spe} that its Cauchy transform satisfies the equation
\begin{equation}\label{eq:semi}
b G(b)-1=\eta\left( G(b)\right)\cdot G(b).
\end{equation}

We want to derive an approximate version of this equation for $G_n$. For this, we will
look at $E[S_n R_n(b)]$.

Let us denote by $S_n^{[i]}$ the version of $S_n$ where the $i$-th variable $X_i$ is
absent, i.e.,
$$S_n^{[i]}:=S_n-\frac 1{\sqrt n} X_i,$$
and by $R_n^{[i]}$ and $G_n^{[i]}$ the corresponding resolvent and Cauchy
transform, respectively, i.e.,
$$R_n^{[i]}(b)=\frac 1{b -S_n\hi}$$
and
$$G_n\hi(b):=E[R_n\hi(b)].$$
For each $i=1,\dots,n$ we have the resolvent identity
\begin{align*}
R_n(b)&=R_n^{[i]}(b)+\frac1{\sqrt n} R_n^{[i]}(b)\cdot  X_i \cdot R_n^{[i]}(b) \\&\quad+
\frac 1n R_n(b)\cdot X_i\cdot R_n^{[i]}(b)\cdot X_i\cdot R\hi_n(b).
\end{align*}

Now we can write
\begin{align*}
E[S_nR_n(b)]&=\sum_{i=1}^n E\bigl[\frac{X_i}
{\sqrt n}\cdot R_n(b)\bigr]\\
&=\sum_{i=1}^n\frac 1{\sqrt n} \Bigl\{ E\bigl[X_i
 \cdot R_n^{[i]}(b)\bigr]
\\&\quad+
\frac 1{\sqrt n} E\bigl[X_i \cdot R_n^{[i]}(b)\cdot X_i \cdot R_n^{[i]}(b)\bigr]\\& \quad
+\frac 1n E\bigl[X_i\cdot R_n(b)\cdot X_i \cdot R_n^{[i]}(b)\cdot X_i\cdot
R\hi_n(b)\bigr]\Bigr\}
\end{align*}

Now we use our assumption that $X_1,X_2,\dots$ are free with respect to $E$, which
implies that $X_i$ is free from $R_n^{[i]}(b)$ with respect to $E$. This implies that
$$
 E[X_i
 \cdot R_n^{[i]}(b)]=
 E[X_i]
 \cdot E[ R_n^{[i]}(b)]=0
$$
and
\begin{align*}
E\left[X_i \cdot R_n^{[i]}(b)\cdot X_i\cdot R_n^{[i]}(b)\right] &= E\left[X_i \cdot E[
R_n^{[i]}(b)]\cdot X_i\right] \cdot E\left[R_n^{[i]}(b) \right]\\&\qquad + E[X_i] \cdot
E\left[R_n^{[i]}(b)\cdot E[X_i]\cdot R_n^{[i]}(b) \right]\\&\qquad
- E[X_i] \cdot E[R_n^{[i]}(b)]\cdot E[X_i]\cdot E[R_n^{[i]}(b)]\\
&= E\left[X_i \cdot E[ R_n^{[i]}(b)]\cdot X_i\right] \cdot E[R_n^{[i]}(b)]\\
&=\eta\left(G\hi_n(b)\right)\cdot G\hi_n(b).
\end{align*}

So we have got finally
\begin{equation}\label{eq:five}
E[S_n R_n(b)]=\frac 1n\left(\sum_{i=1}^n \eta\left(G\hi_n(b)\right)\cdot G\hi_n(b)
+r_1\hi\right),
\end{equation}
where
$$
r_1\hi= \frac 1{\sqrt n} E\left[X_i \cdot R_n(b)
 \cdot X_i\cdot  R_n^{[i]}(b)\cdot
X_i  \cdot R\hi_n(b)\right]
$$

We will now estimate the norm of $r_1\hi$. We could of course just estimate against the
operator norm of $X_i$; however, we prefer, in analogy with the classical case, to do
better without invoking the operator norm and use only as small moments of $X_i$ as
possible.

Note that for our conditional expectation $E$ we have the Cauchy-Schwarz inequality
$$\Vert E[AB]\Vert^2\leq \Vert E[AA^*]\Vert\cdot \Vert E[B^*B]\Vert,$$
and also
$$E[A]^*E[A]\leq E[A^*A]\qquad\text{and}\qquad
E[ABB^*A^*]\leq \Vert BB^*\Vert\cdot E[AA^*]$$ and
$$\Vert E[A]\Vert \leq \Vert A\Vert$$
for any $A,B\in\cA$. Thus, for any $i=1,\dots,n$, we can estimate
\begin{align*}
\Vert E\bigl[X_i\, R_n(b)\,&
 X_i\, R_n^{[i]}(b)\,
X_i\,  R\hi_n(b)\bigr]\Vert^2 \\
&\leq \Vert E\bigl[X_i \, R_n(b)\, R_n(b)^* \, X_i\bigl]\Vert\cdot
\\&\quad\cdot
\bigl\Vert E\bigl[ R_n\hi(b)^*\, X_i \, R_n^{[i]}(b)^* \, X_i\, X_i\, R_n^{[i]}(b)\,
X_i\, R_n\hi(b)\bigl]\bigr\Vert
\end{align*}

We estimate the first factor by
\begin{align*}
\Vert E\bigl[X_i \, R_n(b)\, R_n(b)^* \, X_i\bigl]\Vert &\leq \Vert R_n(b)\Vert^2
\cdot\bigl\Vert
 E\bigl[X_i X_i\bigl]\bigr\Vert\\
&=\Vert R_n(b)\Vert^2\cdot\Vert \eta(1)\Vert\\
&=\alpha_2\Vert R_n(b)\Vert^2
\end{align*}

For the second factor we use again the freeness between $X_i$ and $R_n^{[i]}(b)$. Let us
put
$$R:=R_n^{[i]}(b)$$

Then $X_i$ and $R$ are $*$-free with respect to $E$ and thus, by also invoking
$E[X_i]=0$, we have
\begin{align*}
E[R^*X_iR^*X_iX_iRX_iR]&= E\Bigl[R^*\cdot E\bigl[X_i\, E[R^*]\, X_i\, X_i\, E[R]\,
X_i\bigr]\cdot R\Bigr]
\\
&\quad+
E\Bigl[R^*\cdot\eta\bigl( E[R^*\,\eta(1)\, R] \bigr)\cdot R\Bigr]\\
&\quad- E\Bigl[R^*\cdot \eta\bigl(E[R^*]\,  \eta(1)\, E[R]\bigr)\cdot R\Bigr],
\end{align*}
and thus
\begin{align*}
\left\Vert E\bigl[R^*X_iR^*X_iX_iRX_iR\bigr]\right\Vert&\leq \left\Vert E\Bigl[R^*\cdot
E\bigl[X_i\, E[R^*]\, X_i\, X_i\, E[R]\, X_i\bigr]\cdot R\Bigr]\right\Vert
\\
&\quad+\left\Vert
E\Bigl[R^*\cdot\eta\bigl( E[R^*\,\eta(1)\, R] \bigr)\cdot R\Bigr]\right\Vert \\
&\quad+ \left\Vert E\Bigl[R^*\cdot \eta\bigl(E[R^*]\,  \eta(1)\, E[R]\bigr)\cdot
R\Bigr]\right\Vert
\end{align*}
We estimate
\begin{align*}
&\left\Vert E\Bigl[R^*\cdot E\bigl[X_i\, E[R^*]\, X_i\, X_i\, E[R]\, X_i\bigr]\cdot
R\Bigr]\right\Vert\\ &\qquad\qquad\qquad\qquad\leq\Vert R\Vert \cdot \Vert R^*\Vert\cdot
\bigl\Vert E\bigl[X_i\, E[R^*]\, X_i\, X_i\, E[R]\,
X_i\bigr]\bigr\Vert\\
&\qquad\qquad\qquad\qquad\leq \Vert R\Vert^2\cdot \alpha_4 \cdot \Vert E[R]\Vert\cdot \Vert E[R^*]\Vert\\
&\qquad\qquad\qquad\qquad\leq \alpha_4\cdot \Vert R\Vert^4
\end{align*}
\begin{align*}
\left\Vert E\Bigl[R^*\,\eta\bigl( E[R^*\,\eta(1)\, R] \bigr)\, R\Bigr]\right\Vert \leq
\alpha_2^2\cdot\Vert R\Vert^4,
\end{align*}
and
$$
\left\Vert E\Bigl[R^*\cdot \eta\bigl(E[R^*]\,  \eta(1)\, E[R]\bigr)\cdot
R\Bigr]\right\Vert \leq\alpha_2^2\cdot \Vert R\Vert^4$$

Putting this together yields
$$\left\Vert E\bigl[ R_n\hi(b)^*\, X_i \, R_n^{[i]}(b)^* \, X_i\, X_i\, R_n^{[i]}(b)\,
X_i\, R_n\hi(b)\bigl]\right\Vert \leq (\alpha_4+2\alpha_2^2)\cdot \Vert
R_n^{[i]}(b)\Vert^4,$$ and finally
$$\Vert r_1\hi\Vert \leq \frac 1{\sqrt n}\cdot \sqrt{\alpha_2(\alpha_4+2\alpha_2^2)}\cdot
\Vert R_n(b)\Vert\cdot \Vert R_n^{[i]}(b)\Vert^2.$$

We still need to replace, in \eqref{eq:five}, $G_n^{[i]}(b)=E[R_n^{[i]}(b)]$ by
$G_n(b)=E[ R_n(b)]$. By using the resolvent identity
$$R_n(b)=R_n^{[i]}(b)+\frac 1{\sqrt n} R_n^{[i]}(b)\cdot X_i \cdot R_n(b)$$
we have
$$G_n^{[i]}(b)=G_n(b)+r_2\hi,$$
where
$$r_2\hi:=-
\frac 1{\sqrt n} E[R_n^{[i]}(b)\, X_i \, R_n(b)].$$ As before, we estimate
\begin{align*}
\Vert E[R_n^{[i]}(b)\, X_i \, R_n(b)]\Vert^2 &\leq \Vert E[R_n^{[i]}(b)\, X_i\, X_i\,
R_n^{[i]} (b)^*]\Vert\cdot \Vert E[R_n(b)^* \, R_n(b)]\Vert \\&\leq \alpha_2\cdot \Vert
R_n^{[i]}(b)\Vert^2\cdot \Vert R_n(b)\Vert^2.
\end{align*}

Let us summarize. We have
\begin{align*}
E[S_n R_n(b)]&=\frac 1n \sum_{i=1}^n\left( \eta\left(G\hi_n(b)\right)\cdot G\hi_n(b)
+r_1\hi\right)\\
&=\frac 1n\sum_{i=1}^n\left( \eta\left(G_n(b)+r_2\hi\right)\cdot
\left(G_n(b)+r_2\hi\right) +r_1\hi\right),
\end{align*}
and the estimates
$$\Vert r_1\hi\Vert \leq \frac 1{\sqrt n}\cdot \sqrt{\alpha_2(\alpha_4+2\alpha_2^2)}\cdot
\Vert R_n(b)\Vert\cdot \Vert R_n^{[i]}(b)\Vert^2$$ and
$$\Vert r_2\hi\Vert\leq \frac 1{\sqrt n}\sqrt{\alpha_2}\cdot
  \Vert R_n\hi(b)\Vert\cdot \Vert R_n(b)\Vert.$$
It remains to estimate $\Vert R_n(b)\Vert$ and $\Vert R_n\hi(b)\Vert$. For those we use
the usual estimate for Cauchy transforms (where $\Im b:=(b-b^*)/(2i)$ denotes the
imaginary part of $b$),
$$\Vert R_n(b)\Vert\leq \Vert \frac 1{\Im b}\Vert ,\qquad
\Vert R_n\hi(b)\Vert\leq \Vert\frac 1{\Im b}\Vert.$$ For a formal proof of this estimate,
see, e.g., Lemma 3.1 in \cite{HT}.

We have now
$$
E[S_n R_n(b)]=\eta\left(G_n(b)\right)\cdot G_n(b) +r_3,$$ where
$$r_3=\frac 1n\sum_{i=1}^n\Bigl(\eta(G_n(b))\cdot r_2\hi
+\eta(r_2\hi)\cdot G_n(b)+\eta(r_2\hi)\cdot r_2\hi+r_1\hi\Bigr).$$ Hence
\begin{align*}
\Vert r_3\Vert\leq \frac 1n\sum_{i=1}^n\left( 2\Vert\eta\Vert\cdot \Vert
G_n(b)\Vert\cdot\Vert r_2\hi\Vert + \Vert \eta\Vert\cdot \Vert r_2\hi\Vert^2+\Vert
r_1\hi\Vert\right)\leq{c_n},
\end{align*}
where
$$c_n:=c_n(b):=\frac 1{\sqrt n} \bigl\Vert\frac 1{\Im b}\bigr\Vert^3\sqrt{\alpha_2}
\cdot(2\alpha_2+\sqrt{\alpha_4+2\alpha_2^2})+\frac 1n \bigl\Vert\frac 1{\Im
b}\bigr\Vert^4\alpha_2^2.
$$

Note that $S_nR_n(b)=-1+b R_n(b)$, hence
$$E[S_nR_n(b)]=b G_n(b)-1,$$
and so we finally have found
\begin{equation}\label{eq:G-n}
\eta(G_n(b))\cdot G_n(b)-b G_n(b)+1=-r_3,
\end{equation}
or the inequality:
\begin{equation}\label{eq:ineq}
\Vert \eta(G_n(b))\cdot G_n(b)-b G_n(b)+1\Vert\leq c_n.
\end{equation}

In order to get from this an estimate for the difference between $G_n(b)$ and $G(b)$, we
will now follow the ideas in Section 5 of \cite{HT}, in the improved version from
\cite{HST}.

By \eqref{eq:semi}, we have for all $b\in\cB_+$ the equation
\begin{equation}\label{eq:Gfixed}
b=\frac 1{G(b)}+\eta\bigl(G(b)\bigr)
\end{equation}
for $G(b)$, and, by \eqref{eq:G-n}, the corresponding approximate version for $G_n(b)$:
\begin{equation}\label{eq:Gnfixed}
\Lambda_n(b)=\frac 1{G_n(b)}+\eta\bigl(G_n(b)\bigr),
\end{equation}
where
$$\Lambda_n(b):=b- r_3\cdot G_n(b)^{-1}.$$
A crucial point is now to show that for a sufficiently large set $\tilde O_n\subset
\cB_+$ the quantity $\Im \Lambda_n(b)$ is still positive, so that we can also use
equation \eqref{eq:Gfixed} for $\Lambda_n(b)$. Let us try
\begin{multline*}
\tilde O_n:=\Bigl\{b\in \cB_+\mid c_n(b)<1/2 \quad\text{and}\\ c_n(b)\cdot\Big(\Vert
b\Vert+\alpha_2\cdot\bigl\Vert\frac 1{\Im b}\bigr\Vert\Bigr)\cdot \bigl\Vert\frac 1{\Im
b}\bigr\Vert<1/2 \Bigr\}.
\end{multline*}

The relevance of the condition $c_n(b)<1/2$ is the following: Let us denote
$$B_n(b):=b-\eta(G_n(b)),$$
then inequality \eqref{eq:ineq} takes, for $b\in\tilde O_n$, the form
$$\Vert 1-B_n(b)G_n(b)\Vert \leq c_n(b)<1/2.$$
This, however, implies that $B_n(b)G_n(b)$ is invertible with
$$\Vert G_n(b)^{-1}B_n(b)^{-1}\Vert=\Vert (B_n(b)G_n(b))^{-1}\Vert\leq 2,$$
and thus
\begin{align*}
\Vert G_n(b)^{-1}\Vert&=\Vert G_n(b)^{-1}B_n(b)^{-1}B_n(b)\Vert\\
&\leq 2\Vert B_n(b)\Vert\\
&=2\Vert b-\eta(G_n(b))\Vert\\
&\leq 2\left(\Vert b\Vert +\alpha_2 \cdot\Vert G_n(b)\Vert\right)\\
&\leq 2\left(\Vert b\Vert +\alpha_2\cdot \Vert \frac 1{\Im b}\Vert\right).
\end{align*}

But then the other condition in the definition of $\tilde O_n$ implies that for
$b\in\tilde O_n$ we have
\begin{align}\label{eq:abneu}
\Vert r_3\cdot G_n(b)^{-1}\Vert&\leq \Vert r_3\Vert \cdot\Vert G_n(b)^{-1}\Vert\\&\leq
c_n \cdot 2\left(\Vert b\Vert +\alpha_2 \cdot\Vert \frac 1{\Im b}\Vert\right)<\Vert \frac
1{\Im b}\Vert^{-1}.\notag
\end{align}
Since
$$\Im b\geq \Vert \frac 1{\Im b}\Vert^{-1}\cdot 1,$$
it follows that, for $b\in \tilde O_n$, $\Lambda_n(b)=b-r_3\cdot G_n(b)^{-1}$ is still in
$\cB_+$ and so we can use the equation \eqref{eq:Gfixed} with $\Lambda_n(b)$ as argument,
i.e.,
\begin{equation}\label{eq:GLfixed}
\Lambda_n(b)=\frac 1{G(\Lambda_n(b))}+\eta\bigl(G(\Lambda_n(b))\bigr).
\end{equation}
The point of having both equation \eqref{eq:GLfixed} and equation \eqref{eq:Gnfixed} is
that this implies that
$$G(\Lambda_n(b))=G_n(b).$$
In \cite{HT,HST} this was shown by analytic continuation arguments. We can simplify that
argument by using the fact from \cite{HRS} that the equation
\begin{equation} \label{eq:G}
w=\frac 1G +\eta(G)
\end{equation}
has, for any $w$ with $\Im w>0$, exactly one solution $G\in \cB$ such that $\Im G$ is
negative. Since both $G_n(b)$ and $G(\Lambda_n(b))$ have negative imaginary parts (as
Cauchy transforms at some arguments) and both satisfy the same equation \eqref{eq:G} (for
$w=\Lambda_n(b)$), they must agree.

Then we can, still in the case $b\in \tilde O_n$, estimate in the usual way, by invoking
the resolvent identity:
\begin{align*}
\Vert G_n(b)-G(b)\Vert&=\Vert G(\Lambda_n(b))-G(b)\Vert\\
&=\Vert G(\Lambda_n(b))\cdot (\Lambda_n(b)-b)\cdot G(b)\Vert\\
&\leq \Vert (\Lambda_n(b)-b)\Vert\cdot \Vert G_n(b)\Vert\cdot \Vert G(b)\Vert.
\end{align*}
Both $\Vert G(b)\Vert$ and $\Vert G_n(b)\Vert$ can be estimated by $\Vert 1/\Im b\Vert$
and for the first factor we have, by the second inequality in \eqref{eq:abneu}, that
\begin{align*}
\Vert(\Lambda_n(b)-b)\Vert=\Vert {r_3}{G_n(b)}^{-1}\Vert\leq c_n \cdot 2\left(\Vert
b\Vert +\alpha_2 \cdot\Vert \frac 1{\Im b}\Vert\right)
\end{align*}
Thus, for $b\in \tilde O_n$, we have shown that
\begin{equation}
\Vert G_n(b)-G(b)\Vert\leq c_n \cdot 2\left(\Vert b\Vert +\alpha_2 \cdot\Vert \frac 1{\Im
b}\Vert\right)\cdot \Vert \frac 1{\Im b}\Vert^2
\end{equation}
For $b\in \cB_+\backslash \tilde O_n$, on the other hand, we just use the trivial
estimate
$$\Vert G_n(b)-G(b)\Vert\leq 2\cdot \Vert\frac 1{\Im b}\bigr\Vert$$
together with
\begin{itemize}
\item
if we have $c_n(b)\geq 1/2$, then
\begin{align*}
\Vert\frac 1{\Im b}\Vert &\leq 2c_n\cdot \Vert\frac 1{\Im b}\Vert\\
&\leq 2c_n\cdot \Vert\frac 1{\Im b}\Vert\cdot \Vert b\Vert \cdot \Vert\frac 1{\Im
b}\Vert\\
 &\leq  2c_n\cdot \Vert\frac 1{\Im b}\Vert^2 \cdot \left(\Vert b\Vert +\alpha_2
\cdot\Vert \frac 1{\Im b}\Vert\right)
\end{align*}
\item
if we have $c_n(b)\cdot\left(\Vert b\Vert+\alpha_2\cdot\bigl\Vert\frac 1{\Im
b}\bigr\Vert\right)\cdot \bigl\Vert\frac 1{\Im b}\bigr\Vert\geq 1/2$, then we have again
\begin{align*}
\Vert\frac 1{\Im b}\Vert
 &\leq  2c_n \cdot \left(\Vert b\Vert +\alpha_2
\cdot\Vert \frac 1{\Im b}\Vert\right)\cdot \Vert\frac 1{\Im b}\Vert^2
\end{align*}
\end{itemize}

Thus we have proved the Theorem.

\end{document}